\documentclass[11pt]{article}                      %[20.02.2002]
\makeatletter
\def\draft{\textheight=10.5truein \textwidth=7.5truein \parindent=8pt
           \voffset=-1truein \topmargin=0Truein
           \ifcase \@ptsize \hoffset=-1.5truein \or \hoffset=-1.35truein
                        \or \hoffset=-1.15truein \fi}
\def\quality{\textheight=240mm \textwidth=160mm \topmargin=0Truein
             \ifcase \@ptsize \hoffset=-23mm \or \hoffset=-20mm
                          \or \hoffset=-15mm \fi}
\def\USdraft{\textheight=9.85truein \textwidth=7.5truein \parindent=8pt
             \voffset=-1.0truein \topmargin=0Truein
             \ifcase \@ptsize \hoffset=-1.5truein \or \hoffset=-1.35truein
                          \or \hoffset=-1.15truein \fi}
\def\USquality{\textheight=250mm \textwidth=170mm \topmargin=0Truein
               \voffset=-1truein
               \ifcase \@ptsize \hoffset=-23mm
                       \or \hoffset=-20mm \or \hoffset=-15mm \fi}
\makeatother\USquality

   \newcommand\?[1]{} % comments
\def\beq#1#2{\begin{equation} \label{#1} #2 \end{equation}}
\def\bea#1{\begin{eqnarray*} #1 \end{eqnarray*}} \def\a{\!\!\!&\!\!\!\!&}
\def\function#1{\left\{\!\!\!\begin{array}{ll} #1 \end{array} \right.}
\def\n{\noindent} \def\CR{$$ $$} \def\dist{\, {\rm dist}}
\def\proof{\smallskip \noindent {\bf Proof. \ }}       %start of proof
\def\blanksquare{\,\,\,$\sqcup\!\!\!\!\sqcap$}         %blank  square
\def\qed{\hfill\blanksquare\linebreak\smallskip\par}   %end of proof

\def\bline(#1,#2)(#3,#4)(#5){\put(#1,#2){\line(#3,#4){#5}}}  %straight line

\def\thname{Theorem}     \def\lmname{Lemma}      \def\prname{Proposition}
\def\dfname{Definition}  \def\crname{Corollary}  \def\rmname{Remark}

\newtheorem{theorem}{\thname}[section]   %Numbering: Theorem--Other section
\newtheorem{lemma}{\lmname}[section]     %{lemma}[theorem]{Lemma}   section
\newtheorem{corollary}[lemma]{\crname}   %lemma

\newtheorem{dftn}{\dfname}[section]
\newtheorem{rmrk}[lemma]{\rmname}
\newenvironment{remark}{\begin{rmrk}\rm}{\end{rmrk}}     %lemma

             \catcode`@=11 \@addtoreset{equation}{section} \catcode`@=12
\newcommand\mlbscale{1pt} %to change: \renewcommand\mlbscale{1.3pt}
\newif\iffigs\figstrue %\newif\iffigs\figsfalse -- to fake figures
\def\Bfig(#1,#2)#3#4{\begin{figure} \begin{center}
    \setlength{\unitlength}{\mlbscale}
       \iffigs \begin{picture}(#1,#2) #3 \end{picture}
       \else \begin{picture}(60,10)(0,0)
                   \put(0,0){\framebox(60,10){Figure}} \end{picture} \fi
    \end{center} \caption{#4} \end{figure}}
\def\bpic(#1,#2)#3{\setlength{\unitlength}{\mlbscale}
    \begin{picture}(#1,#2) #3 \end{picture}}
\def\IR{\hbox{\rm I\kern-.2em\hbox{\rm R}}} \def\ep{\varepsilon}
\def\IZ{\hbox{{\rm Z}\kern-.3em{\rm Z}}}    \def\const{\, {\rm Const} \,}

\def\map{T}

\def\cM{{\cal M}}

\def\noprint#1{}

\def\SBR{SRB\ } \def\cM{{\cal M}} \def\cN{{\cal N}} \def\cL{{\cal L}}
\def\cmap{{\cal\map}}  \def\cI{{\cal I}} %interaction
\def\toas#1{\stackrel{#1}{\longrightarrow}}   

\begin{document}
%%%%%%%%%
\title{Multicomponent dynamical systems:\\ \SBR measures and phase
transitions}

\author{Michael Blank\thanks{Russian Academy of Sciences,
                             Inst. for Information Transmission Problems,
                             and Georgia Institute of Technology;
                             e-mail: blank@obs-nice.fr},
        Leonid Bunimovich\thanks{School of Mathematics,
                                 Georgia Institute of Technology;
                                 e-mail: bunimovh@math.gatech.edu}}
\date{February 20, 2002}
\maketitle

\n{\bf Abstract.} We discuss a notion of phase transitions in
multicomponent systems and clarify relations between
deterministic chaotic and stochastic models of this type of
systems. Connections between various definitions of \SBR measures
are considered as well.

\section{Introduction}

The aim of the present paper is twofold: to study the notion of
phase transitions in multicomponent systems and to clarify the
relations between deterministic chaotic and stochastic models of
this type of systems. We also discuss the differences in the
approaches to multicomponent systems in statistical physics and
in dynamical systems theory. In the former case it is basically a
system of interacting particles, while in the later case each
component of a multicomponent systems can have a nontrivial
(local) dynamics, which leads to more rich evolution and thus to
more rich statistical properties. Our definition of a
multicomponent system (see Section 4) assumes that it consists of
components (local susbsystems) which have their intrinsic
dynamics and besides interact with each other. In particular, we
argue that in distinction to situations considered in statistical
physics, phenomenona similar to phase transitions can appear even
in systems with a finite number of components (degrees of
freedom). Indeed, in the absence of interactions dynamics of
particle systems is trivial. Besides, some interesting phenomena,
even in the presence of interactions, may appear only in the
limit when the number of particles goes to infinity. Moreover,
often additional assumptions like indecomposability (see e.g.
\cite{MacK}), are used to emphasize that the situation under
study is impossible in a finite dimensional setting. On the other
hand, as we shall show in Section~\ref{s-SBR}, already the
simplest one-dimensional dynamical systems satisfying the
indecomposability assumption (and even the assumption of
topological transitivity) may be non ergodic, which shows that
restrictions of this type are not quite reasonable in the context
of general dynamical systems.

A recent progress in the analysis of chaotic spatially extended
dynamical systems allows to advance in answering to a long
standing question how to define exactly phase transitions
rigorously and what are the conditions for phase transitions in
this type of systems and more generally in multicomponent
systems. The problem with the definition of the phase transition
phenomenon is that this notion is used in different ways in
statistical physics (see review in \cite{EFS}), moreover various
existing approaches to this notion for spatially extended systems
\cite{Bl21,Bl-mon,Bun-C,MacK} lead to different statements about
the existence of phase transitions. As it was already mentioned
in \cite{Bun-C,Bl21,MacK} it is essential to make distinction
between qualitative changes in the topological behaviour of a
system, called bifurcations in the dynamical systems theory, and
changes to measure-theoretical properties which we shall identify
with phase transitions. In the present paper by the phase
transition we shall mean a change of a number of \SBR measures
(see the definition and discussion further) which one can
naturally identify with phases in statistical physics. In fact,
this is the same general idea, which was used earlier in a number
of papers \cite{Bun-C,Bl21,MacK}, however differences in the
definitions of the \SBR measures (which we shall discuss in
detail) lead again to different statements about the phase
transitions.

If local components of our multicomponent system are identical
and the interaction is translationally invariant one can discuss
finite dimensional approximations $(T^{(d)},X^d)$ (where $X^{d}$
is a direct product of $d$ identical copies of $X$) with certain
boundary conditions (e.g. with periodic boundary conditions).
Assuming that we are able to study ergodic properties of those
finite dimensional approximations for any $d<\infty$, one of the
major problems is to analyze how their limiting (as $d\to\infty$)
behavior corresponds to the dynamics of the entire (infinite)
system. In particular, it might be possible that for each finite
$d$ there is only one \SBR measure $\mu_d$, but their limit
points do not coincide with the \SBR measures of the
multicomponent system.

\section{\SBR measures}\label{s-SBR}

Let $(X,\rho)$ be a compact metric space with a certain reference
measure $m$ on it, e.g. a finite dimensional unit cube or torus
with the Lebesgue measure. Consider a nonsingular with respect to
the measure $m$ map $\map$ from the space $X$ into itself (i.e.
$m(\map^{-1}A)=0$ whenever $m(A)=0$). The pair $(\map,X)$ defines
a deterministic dynamical system. To study statistical properties
of this dynamical system we need to consider its action in the
space of measures. Let $\cM(X)$ be the space of probabilistic
measures on $X$ equipped with the topology of weak convergence of
measures. Then the induced map $\map^*:\cM(X)\to\cM(X)$ is
defined as follows:%
\beq{def:ind-map}{\map^*\mu(A):=\mu(\map^{-1}A)}%
for any measure $\mu\in\cM(X)$ and any Borel set $A\subseteq X$.

We shall call a measure $\mu_\map$ a {\em natural} measure for
the map $\map$ if there exists an open subset $U\subseteq X$
(called the basin of attraction for the measure $\mu_\map$) such
that for any measure $\mu\in\cM(X)$ absolutely continuous with
respect to the reference measure $m$ and having its support in
$U$ we have:%
\beq{def-nat}{\frac1n \sum_{k=0}^{n-1}{\map^*}^k\mu
              \toas{n\to\infty} \mu_\map}%
In other words, the measure $\mu_\map$ is a stable fixed point of
the dynamics of absolutely continuous initial measures. A similar
definition has been used e.g. in \cite{FOY,Bl-mon}.

Observe that from the point of view of the action of the map
$\map$ in the space of measures a natural measure (not necessary
unique) is nothing more than a stable (with respect to absolutely
continuous initial conditions in the space of measures) fixed
point of the induced map $\map^*$, i.e. an attractor. This object
is well known in ergodic theory of dynamical systems and
corresponds to one of the definitions of \SBR measures (see e.g.
\cite{Bl-mon,Bun-S,FOY,MacK}), which we shall discuss in a
moment. Observe also that one of the advantages of this
definition is that it works without any changes for true random
Markov chains as well. Indeed, let $\map^*$ be the transfer
operator of a Markov chain with the phase space $X$. This
operator generates the dynamics of measures on $\cM(X)$, i.e.
this is a conjugate operator to the Markov operator (transition
matrix) of the Markov chain under consideration. Then the
relation~(\ref{def-nat}) defines the notion of the natural
(\SBR\hskip-1.5mm) measure in this true random setting as well.

In the literature one can find three main different approaches to
a definition of \SBR measures in the deterministic setting. We
already mentioned one of them (natural measure).

The second definition is very close to the previous one, with the
only difference that one considers pathwise convergence of sample
measures, rather than a convergence of orbits of the induced map.
Namely, by \SBR measure in this case one means the common limit as
$n\to\infty$ for $m$-almost all points $x\in U$ of
$\frac1n\sum_{k=1}^{n-1}\delta_{\map^n x}$, where $\delta_x$
stands for the $\delta$-measure at the point $x$.

Observe that both these definitions are based on Gibbs idea of
construction of stationary measures and therefore they both are
closely related to the statistical physics formalism and to the
well known Bogolyubov-Krylov approach in dynamical systems theory.

The third definition of \SBR measure is based on a completely different
observation, namely that for some `good' dynamical systems with strong
stochastic properties (for example, uniformly hyperbolic systems) an
\SBR measure, which corresponds to any of the above definitions, has a
marginal
distribution absolutely continuous with respect to the Lebesgue measure on
the
so called unstable foliation of the dynamical system.
This is indeed a very important statistical feature of the dynamics, but it
is well defined only in the case of hyperbolic dynamical systems and
therefore
it is not clear what is the reason to use it as a definition of the \SBR
measure for more general dynamical systems. Denote by $\mu_\map$,
$\mu_\map^p$
and $\mu_\map^c$ the versions of \SBR measures corresponding to the three
definitions given above respectively. The following simple statement
describes
connections between these definitions. It shows also that the first version,
which we call the natural invariant measure is, indeed, more natural than
the
others.

\begin{theorem}\label{t:sbr}%
\begin{itemize}
\item[(a)] If $\mu_\map^p$ exists then $\mu_\map$ is also well
    defined and $\mu_\map=\mu_\map^p$.%
\item[(b)] If $\mu_\map$ is ergodic and invariant then
    $\mu_{\map}^{p}$ is well defined and $\mu_\map^p=\mu_\map$,
    however without ergodicity this might not hold.%
\item[(c)] If $\mu_\map^c$ is stable with respect to the dynamics of
    absolutely continuous initial measures then $\mu_\map=\mu_\map^c$.%
\item[(d)] The existence of the $\mu_\map$ or $\mu_\map^p$ does not
    imply the existence of $\mu_\map^c$.
\end{itemize}
\end{theorem}

\proof We start from the assertion (a). By definition for a $m$-a.a.
point $y\in U\subseteq X$ we have the week convergence of the sequence of
measures $\frac1n\sum_{k=1}^{n-1}\delta_{\map^ky}\to\mu_\map^p$, i.e. for
any
continuous function $\phi:X\to\IR^1$ and $m$-a.a. $y\in U$ we have
$$ \int \phi(x)~d\left(\frac1n\sum_{k=1}^{n-1}\delta_{\map^ky}\right)
  = \frac1n\sum_{k=1}^{n-1}\phi(\map^ky) \toas{n\to\infty} \int
\phi~d\mu_\map^p .$$
Choose a measure $\mu\in\cM(X)$ with a support in $U$ being
absolutely continuous with respect to $m$ and consider its Cesaro
averages: $\mu_n:=\frac1n\sum_{k=1}^{n-1}{\map^*}^k\mu$. Then
using the above convergence and absolute continuity of the
measure $\mu$ we get%
\bea{ \int \phi~d\mu_n
  \a = \int \phi ~d\left(\frac1n\sum_{k=1}^{n-1}{\map^*}^k\mu\right)
     = \int \frac1n\sum_{k=1}^{n-1}\phi(\map^kx) ~d\mu \\
  \a \toas{n\to\infty} \int (\int \phi~d\mu_\map^p)~d\mu = \int
\phi~d\mu_\map^p ,}%
which proves the first assertion.

It is of interest that even if the measure $\mu_\map^p$ exists and is unique
it might be non ergodic. To show this consider the following example,
proposed
by G.~Del Magno:
\beq{ex:GiGi}{
       \map x := \function{(1-\sin(\pi x-\pi/2))/2  & \mbox{if } 0< x < 1,
\\
                           x & \mbox{if } x\in\{0,1\} } .}
One can easily show that the locally maximal attractor in this
example consists of two fixed points at 0 and 1 and that for any
initial point $x\in(0,1)$
$$ \frac1n\sum_{k=1}^{n-1}\delta_{\map^n x}
   \toas{n\to\infty}\frac12(\delta_0+\delta_1) = \mu_\map^p .$$
On the other hand, this measure is nonergodic, since both points
0 and 1 are fixed points.

The first part of the assertion (b) follows immediately from
Birkhoff ergodic theorem. Observe that in order to apply this
theorem we need the measure $\mu_{\map}$ to be invariant. If the
map $\map$ is continuous this is certainly correct (it is enough
to apply the induced map to the both sides of the limit
construction in the definition of $\mu_{\map}$), but for a
general nonsingular map this fact is not an immediate consequence
of the definition of the natural measure. To demonstrate this
consider the following one dimensional map from the unit interval
into itself: $\map x=x/2$ for all $x\in(0,1]$ with the only
discontinuity at the origin: $\map0=1$. In this example under the
action of the induced map $\map^{*}$ any probabilistic measure
converges to the $\delta$-measure at $0$, but the map $\map$ has
no invariant measure at all.

To finish the proof of this assertion we need to show that the
natural measure might not coincide with $\mu_\map^p$, which at
first sight looks rather doubtful. Observe, that the example of
the map (\ref{ex:GiGi}) shows that even being unique the measure
$\mu_\map$ might be nonergodic. Indeed, again by the same
argument as above the images of any absolutely continuous
probabilistic measure converge to $\frac12(\delta_0+\delta_1)$,
which is nonergodic. Still in this case the measures $\mu_\map$
and $\mu_\map^p$ coincide.

To demonstrate that these two measures might not coincide, we
consider the following one-dimensional map introduced recently in
\cite{Inoue}: \beq{ex:inoue}
    {\map x := \function{x+4x^3     & \mbox{if } 0  \le x < 1/2  \\
                         x-4(1-x)^3 & \mbox{if } 1/2\le x \le1 } .}
This map has two neutrally unstable fixed points $0$ and $1$. It
has been shown \cite{Inoue} that for any sufficiently small
$\delta>0$ for Lebesgue a.a. points $x\in[0,1]$:
$$ \limsup_{n\to\infty} \frac1n\sum_{k=0}^{n-1}1_{[0,\delta]}(\map^k x)=1,
\qquad
   \liminf_{n\to\infty} \frac1n\sum_{k=0}^{n-1}1_{[0,\delta]}(\map^k x)=0
,$$
where $1_A$ stands for the characteristic function of the set
$A$. Therefore a measure $\mu_\map^p$ does not exist in this
case. On the other hand, in \cite{Kel} it has been shown that for
any probabilistic absolutely continuous measure $\mu$ the
sequence of measures ${\map^*}^n\mu \to (\delta_0 +
\delta_1)/2=:\mu_T$ weakly as $n\to\infty$, which proves that the
natural measure $\mu_T$ is well defined. Observe, however, that
there are invariant sets $\{0\}$ and $\{1\}$ having
$\mu_T$-measure 1/2 each, which contradicts to ergodicity as in
the example (\ref{ex:GiGi}).

The statement (c) is a trivial corollary to the definition of the natural
measure. The last assertion (d) follows from the observation that in the
case
of a globally contracting map the measures $\mu_\map$ and $\mu_\map^p$
coincide with a $\delta$-measure being the only invariant measure of the
system, while the measure $\mu_\map^c$ does not exist. \qed

Observe that the measure $\mu^c_\map$ may also exist but be not
finite. It happens already for the most closest to uniformly
hyperbolic so called almost Anosov diffeomorphisms, i.e. for
diffeomorphisms which are uniformly hyperbolic away from a finite
set of points \cite{HY}. Similar results for the case of
one-dimansional neutral maps are also well known (see e.g.
\cite{Inoue} and references therein).

In the literature dedicated to phase transitions one can find the
assumption of indecomposability (see e.g. \cite{MacK}) introduced
to stress the necessity to consider only infinite dimensional
systems. Roughly speaking the indecomposability means the
following: if we have two finite pieces of trajectories
$\{\map^kx\}_{k=1}^n$ and $\{\map^ky\}_{k=1}^m$ on the locally
maximal attractor of the map $\map$, then $\forall\ep>0$ there
exists a point $z=z(\ep)\in X$ and an integer $N$ such that
$\rho(\map^kx,\map^kz)<\ep$ for any $k=1,\dots,n$ and
$\rho(\map^ky,\map^{N+k}z)<\ep$ for any $k=1,\dots,m$. It is not
hard to understand that if the map $\map$ is continuous this is
equivalent to the assumption of {\em topological transitivity} of
the map, i.e. that for any two open subsets $A,B$ having nonempty
intersections with the locally maximal attractor there exists a
number $N$ such that $\map^N A\cap B\ne\emptyset$, which in turn
is equivalent to the presence of a trajectory densely covering
the locally maximal attractor.
%On the other hand, if the map is not continuous those properties
%fail to be equivalent and, in fact, the assumption of topological
%transitivity looks more reasonable.

Let us show now that in distinction to models of classical
statistical physics even elementary one-dimensional dynamical
system might be nonergodic but topologically transitive. We
formulate this statement as a lemma, but in fact it is already
proven above.

\begin{lemma}\label{l:transitive} The maps (\ref{ex:GiGi}) and
(\ref{ex:inoue}) are topologically transitive but their unique
\SBR measures are non ergodic.
\end{lemma}

\section{Deterministic models of Markov chains}\label{s:det-mod}

The aim of this section is to study connections between Markov
chains and piecewise linear maps.

Let $\map$ be a nonsingular map from the unit interval $[0,1]$
into itself and let $\bar\Delta:=\{\Delta_i\}_{i=1}^W$ be a
partition of $[0,1]$ into disjoint intervals. The number $W$ of
elements of the partition $\Delta$ might be infinite. This
partition is called a {\em special} partition for the map $\map$
if the restriction of the map $\map$ to the inner part of any
interval $\bar\Delta_i$ is a diffeomorphism onto its image. The
map $\map$ is called {\em Markov} if there exists a special
partition $\Delta$ (which is also called Markov one) for which
the following property holds: $\map\Delta_i = \cup_{j\in
I(i)}\Delta_j$ for each $i$. Given a Markov partition
$\bar\Delta$ consider a subset of a set of probabilistic measures
$$ \cM_u([0,1], \bar\Delta):=\{\mu\in\cM([0,1]):\quad
   \frac{\mu(I)}{\mu(I')} = \frac{|I|}{|I'|}, \quad
   \forall i ~{\rm and}~ \forall ~{\rm intervals}~ I,I'\subset\Delta_i\} ,$$
where $|I|$ stands for the Lebesgue measure of the interval $I$.
In other words, $\cM_u([0,1], \bar\Delta)$ corresponds to the set
of piecewise uniform distributions on intervals $\bar\Delta$, i.e.
the restriction of any measure from this set to an interval
$\Delta_i$ for each $i$ is proportional to the Lebesgue measure.
Observe that if the map $\map$ is Markov and is mixing, then the
natural measure is unique and its density with respect to the
reference measure is a piecewise constant function on the elements
of the partition.

\begin{theorem}\label{t:det-mod1} For any given transition matrix $P$
of a Markov chain with $N$ states (where $N$ can be infinite)
there exists a piecewise linear one-dimensional map $\map$ with a
Markov partition $\bar\Delta'$ such that the restriction of the
induced map $\map^*$ to $\cM_u([0,1], \bar\Delta')$ is equivalent
to the left action of the matrix $P$ in the space of
distributions.
\end{theorem}

\proof Let $\bar\Delta$ be any partition of the unit intervals
into $N$ subintervals. (For example, if $N<\infty$ one can choose
a partition into $N$ equal intervals, while if $N=\infty$ one can
consider a countable partition into intervals
$\Delta_i:=(2^{-i},2^{-i+1}]$ for $i=1,2,\dots$.) For a given
integer $i$ let $I_i:=\{i_j\}_j$ be the collection of indices
such that $P_{ii_j}>0$. Consider now a subpartition of the
interval $\Delta_i$ into intervals $\{\Delta_{ij}\}_j$ of lengths
$|\Delta_{ij}| = P_{ii_j}|\Delta_{i}|$. This refined partition,
consisting of intervals $\Delta_{ij}$, will be a special
partition $\Delta'$ for a piecewise linear map
$\map:[0,1]\to[0,1]$ defined as follows: on each interval
$\Delta_{ij}$ the map $\map$ is a linear map from this interval
onto the interval $\Delta_{i_j}$. On the other hand, since for
any pair of indices $i,j$ the image under the action of the map
$\map$ of the interval $\Delta_{ij}$ is an interval from the
original partition $\bar\Delta$, and hence a union of intervals
from the partition $\bar\Delta'$, we indeed have shown that this
partition is Markov.

It remains to discuss the equivalence of the restricted induced
map $\map^{*}$ with the left action of the transition matrix $P$.
Observe that the density of a measure $\mu \in \cM_u([0,1],
\bar\Delta')$ is well defined and is a piecewise constant
function on intervals $\Delta_{i}'$. Associate to the measure
$\mu$ a vector $\bar p(\mu)$ with components $(\bar p(\mu)))_{i}
:= \mu(\Delta_{i}')$. Since the map $\map$ constructed above is
Markov, each of the intervals $\Delta_{i}'$ is mapped under the
action of $\map$ onto a union of intervals from the partition
$\bar\Delta'$. Moreover, the map $\mu$ is linear on each element
of the partition and thus $(\bar p(\map^{*}\mu))_{i}=(\bar
p(\mu)P)_{i}$ for each index $i$. \qed

\begin{corollary} Let a number of states $N$ of the Markov chain be finite
and let for each $i$ a number of positive elements in the $i$-th
row of the matrix $P$ do not exceed $K$. Then the number of
elements in the special partition of the one-dimensional map
constructed above is at most $NK$.
\end{corollary}

One might ask if it is possible to construct a continuous version
of the map representing say a finite Markov chain. Indeed, at the
first sight, it looks like using a rearrangement of the elements
of the partition $\Delta$ and changing their lengths this should
be possible. However the following example of a 3-state Markov
chain shows that in general there is no continuous deterministic
model of this type and illustrates also the procedure described
in the proof of above result.

Consider a Markov chain with three states and the following transition
matrix:%
\beq{tr-matr}{
   P = \left(\matrix{ \frac12 & \frac12 & 0 \cr
                      0       & \frac12 & \frac12 \cr
                      \frac12 & 0       & \frac12 \cr }\right) .}%
Then the construction described in the proof of
Theorem~\ref{t:det-mod1} leads to the following map (see
Fig.\ref{model1} below):
$$ \map x := \function{2x     & \mbox{if } 0  \le x<1/3  \\
                       2x-1/3 & \mbox{if } 1/3\le x<2/3  \\
                       2x-4/3 & \mbox{if } 2/3\le x<5/6  \\
                       2x-1   & \mbox{if } 5/6\le x\le1 .}$$
%%%%%%%%%%%%%%%%%%%%%%%%%%%%%%%%%%%%%%%%%%%%%%%%%%%%%%%%%%%%%%%%%%
%% Counterexample to continuity of the model
\Bfig(150,150)
      {\bline(0,0)(1,0)(150)   \bline(0,0)(0,1)(150)
       \bline(0,150)(1,0)(150) \bline(150,0)(0,1)(150)
       %\bline(0,0)(1,1)(150)
       \bezier{50}(0,50)(75,50)(150,50)
       \bezier{50}(0,100)(75,100)(150,100)
       \bezier{50}(50,0)(50,75)(50,150)
       \bezier{50}(100,0)(100,75)(100,150)
       \bezier{50}(125,0)(125,75)(125,150)
       %\bline(0,50)(1,0)(150)  \bline(0,100)(1,0)(150)
       %\bline(50,0)(0,1)(150)  \bline(100,0)(0,1)(150)
       \thicklines
       \bline(0,0)(1,2)(50)   \bline(50,50)(1,2)(50)
       \bline(100,0)(1,2)(25)  \bline(125,100)(1,2)(25)
       \thinlines
       \put(-25,20){}  \put(-25,70){$\Delta_{2}$}
       \put(-25,110){$\Delta_{3}$} \put(-25,130){$\Delta_{4}$}
       \put(20,-10){$\Delta_{1}$}  \put(70,-10){$\Delta_{2}$}
       \put(110,-10){$\Delta_{3}$} \put(130,-10){$\Delta_{4}$}
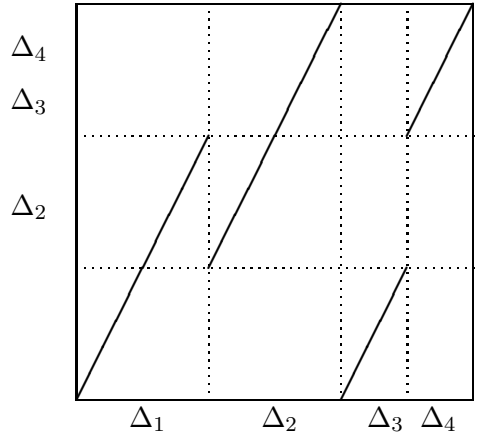
      }{Deterministic model of the Markov
chain~(\ref{tr-matr}).\label{model1}}
%%%%%%%%%%%%%%%%%%%%%%%%%%%%%%%%%%%%%%%%%%%%%%%%%%%%%%%%%%%%%%%%%%
In this case we can make a smaller number of the elements of the
partition, i.e. 4 instead of 6, but a simple geometric argument
shows that one cannot rearrange these intervals to make the map
be continuous, even if we consider it on a circle instead of the
interval. Indeed, the structure of the transition matrix
(\ref{tr-matr}) $P$ is such that the interval(s) corresponding to
at least one of the tree states of the Markov chain should be
mapped to the intervals corresponding to the two other states,
inevitably having a gap between them.

Another and a more interesting example is the deterministic model
of a random walk on nonnegative integers, defined for any
positive $i$ by transition probabilities $p_{i,i-1}, p_{i,i},
p_{i,i+1}$ to go to the left, to remain in the current position,
and to go to the right respectively, and for $i=0$ by transition
probabilities $p_{0,0}$ and $p_{0,1}$. The corresponding map is
shown on Fig.\ref{model2}(a).
%%%%%%%%%%%%%%%%%%%%%%%%%%%%%%%%%%%%%%%%%%%%%%%%%%%%%%%%%%%%%%%%%%
%% Random walk
\Bfig(350,150)
      {\bline(0,0)(1,0)(150)   \bline(0,0)(0,1)(150)
       \bline(0,150)(1,0)(150) \bline(150,0)(0,1)(150)
       \bezier{75}(0,0)(75,75)(150,150)
       \bezier{50}(0,40)(75,40)(150,40)
       \bezier{50}(0,60)(75,60)(150,60)
       \bezier{50}(0,83)(75,83)(150,83)
       \bezier{50}(0,98)(75,98)(150,98)
       \bezier{50}(40,0)(40,75)(40,150)
       \bezier{50}(60,0)(60,75)(60,150)
       \bezier{50}(83,0)(83,75)(83,150)
       \bezier{50}(98,0)(98,75)(98,150)
       \bline(60,83)(-1,-3)(20)
       \bline(60,40)(2,5)(23)
       \bline(83,60)(1,4)(15)
       \bline(0,0)(1,3)(7)  \bline(150,150)(-1,-3)(7)
       \put(15,18){$\dots$} \put(120,125){$\dots$}
       \put(40,-10){$\Delta_{i+1}$}
       \put(67,-10){$\Delta_{i}$}
       \put(83,-10){$\Delta_{i-1}$}
       \put(-25,45){$\Delta_{i+1}$} \put(-25,68){$\Delta_{i}$}
       \put(-25,90){$\Delta_{i-1}$}
       \put(70,-25){\bf(a)}
      %standard random walk
       \put(200,0){\bpic(120,150){
       \bezier{50}(10,-5)(10,75)(10,157)
       \bezier{50}(40,-5)(40,75)(40,157)
       \bezier{50}(70,-5)(70,75)(70,157)
       \bezier{50}(100,-5)(100,75)(100,157)
       \bezier{50}(0, 0)(75, 0)(110, 0)
       \bezier{50}(0,40)(75,40)(110,40)
       \bezier{50}(0,90)(75,90)(110,90)
       \bezier{50}(0,150)(75,150)(110,150)
       \bline(10, 0)(3,4)(30) \bline(40,40)(3,5)(30) \bline(70,90)(1,2)(30)
       \put(10,-10){$\Delta_{i,i+1}$} \put(45,-10){$\Delta_{i,i}$}
       \put(70,-10){$\Delta_{i,i-1}$}
       \put(-20,20){$\Delta_{i+1}$} \put(-20,65){$\Delta_{i}$}
       \put(-20,120){$\Delta_{i-1}$} \put(50,-25){\bf(b)}
                  }}
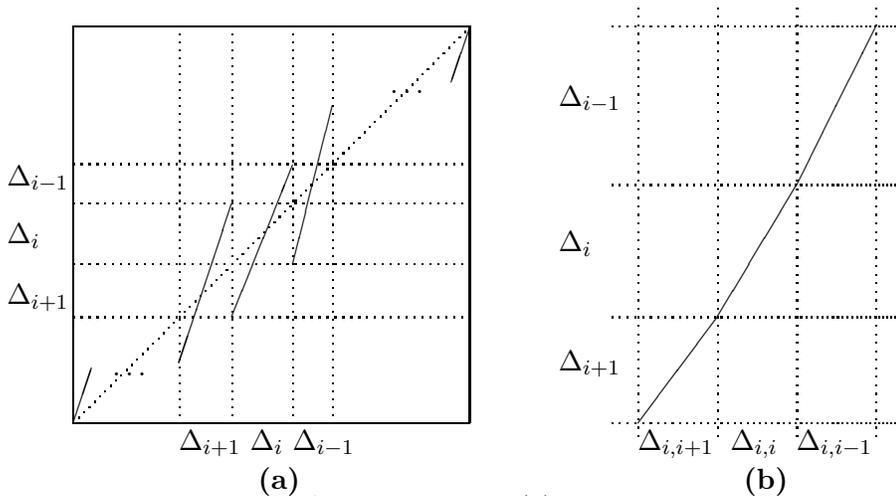
      }{Deterministic models of random walks. (a) nonhomogeneous random
walk,
        (b) homogeneous random walk with
            $\Delta_{i}=\Delta_{i,i-1}\cup\Delta_{i,i}\cup\Delta_{i,i+1}$.
        \label{model2}}
%%%%%%%%%%%%%%%%%%%%%%%%%%%%%%%%%%%%%%%%%%%%%%%%%%%%%%%%%%%%%%%%%%
Observe that $p_{i,j}=|\Delta_{j}|/|\map\Delta_{i}|$ for any pair
of indices $i,j$ such that $|i-j|\le1$.

It is not hard to realize also a (space) homogeneous random walk
on nonnegative integers which is define by transition
probabilities $p_{i,i-1}=p_{L}, p_{i,i}=1-p_{L}-p_{R},
p_{i,i+1}=p_{R}$ for all $i>0$ and $p_{0,0}=1-p_{R},
p_{0,1}=p_{R}$. The corresponding map is shown on
Fig.\ref{model2}(b) and the transition probabilities are equal to
$p_{L}=|\Delta_{i,i-1}|/|\Delta_{i}|$ and
$p_{R}=|\Delta_{i,i+1}|/|\Delta_{i}|$, where
$\Delta_{i}=\Delta_{i,i-1}\cup\Delta_{i,i}\cup\Delta_{i,i+1}$.
Observe that in this case the restriction of the map to
$\Delta_{i}$ is only piecewise linear (in distinction to the
previous case).

While being sufficiently general Markov chains with a countable
number of states do not describe the dynamics of the so called
probabilistic cellular automata, which we are going to discuss
now. Let $G$ be a finite or countable graph, and let $v$ be a
function with a finite number of values (not larger than
$K<\infty$) defined on the vertices of this graph, which we shall
denote by the same letter $G$, i.e. $v:G\to\{1,2,\dots,K\}$. We
assume also that the graph $G$ is locally finite, i.e. for a
vertex $g\in G$ we define its neighborhood $O(g)$ as the union of
vertices to which it is connected. We assume that the graph $G$ is
locally finite, i.e. for each vertex $g\in G$ its neighborhood
$O(g)$ consists of at most $L<\infty$ vertices of the graph $G$.
A probabilistic cellular automaton on the graph $G$ is defined as
follows. For each $g\in G$, each value of $v(g)$, and each
configuration of values of the function $v(O(g))$ we define
transition probabilities $p(\cdot,i)$ for $i\in\{1,2,\dots,K\}$
and thus the dynamics of values $v(g)$. (Observe that the number
of these configurations is finite).

Let $P_M$ be the transfer operator corresponding to a Markov
chain $M$ acting on the state space $X_M$, i.e.
$P_M:\cM(X_M)\to\cM(X_M)$. We shall say that this Markov chain is
{\em equivalent} to a dynamical system $(\map, X)$ if there is a
subspace $\cM_M(X)$ of the space of probabilistic measures on $X$
and a homeomorphism $\pi:\cM(X_M)\to\cM_M(X)$ such that for any
two measures $\mu\in\cM(X_M)$ and $\nu\in\cM_M(X)$ we have
$\pi(P_M\mu)=\map^*(\pi\mu)$ and $\pi P_M(\nu)=\map^*(\pi\nu)$.
Observe that the equivalence of the restricted induced map to the
left action of the transition matrix in Theorem~\ref{t:det-mod1}
is a special case of this general definition.

\begin{theorem}\label{t:mod-det2} For any probabilistic cellular
automaton on a locally finite graph $G$ there exists a
deterministic dynamical system described by a countable number of
maps, whose dynamics is equivalent to the dynamics of the
probabilistic cellular automaton.
\end{theorem}

\proof The idea of the construction of the equivalent dynamical
system for a probabilistic cellular automaton is the following:
in each vertex of the graph we define a number of maps
(corresponding to all possible configurations of states of
neighbors and itself). A choice of the map is given by the
configuration of these states.

Observe, that for a given vertex $g\in G$ a total number of
various configurations of values of the function $v$ in its
neighborhood $O(g)$ cannot exceed $K^L<\infty$. On the other
hand, for any given configuration of $v(O(g))$ the construction
in the proof of the previous theorem can be applied to define a
piecewise linear map $\map_{v(O(g))}$ from a unit interval into
itself with the Markov partition consisting of at most $K^K$
intervals $\Delta_i$, which has the dynamics ``equivalent'' to
the transition probabilities of our cellular automaton. Consider
now a multicomponent system whose phase space $X:=[0,1]^G$ is a
direct product of unit intervals, at each vertex $g\in G$ we have
a finite number of one-dimensional maps $\map_{v(g)}$ and for any
point $\bar x=(x_g)\in X$ we define its image as follows. Let
$x_{g'}\in\Delta_{i(g')}$ for any $g'\in O(g)$ then the map
$\map_{v(O(g))}(x_g)$, corresponding to the configuration
$\{i(g')\}$, defines the value of the $g$-th coordinate of $\bar
x$. \qed

\begin{remark} Using the same argument as above we can construct
a deterministic model with $K=\infty$ if we assume additionally
that for any $i\in\{1,2,\dots,K\}$ only a finite number of
transition probabilities (in the definition of the cellular
automaton) are positive, however the assumption of the local
finiteness of the graph $G$ we cannot drop.
\end{remark}

A particular case of this construction has been used in
\cite{MacK}. Note that our description of the deterministic model
corresponding to an arbitrary probabilistic cellular automaton is
somewhat simpler and more explicit than the one in proposed in
the cited paper.

Let us consider briefly an inverse problem -- an approximation of
a general dynamical system by a finite (countable) state Markov
chain. In distinction to the case of piecewise linear Markov maps
(considered above) in the general case one cannot construct an
``equivalent'' finite or countable state Markov chain, however an
approximation is still possible. One of the simplest approaches
here was proposed in \cite{Ul}. Let $\map$ be a nonsingular map
from a compact metric space $(X,\rho)$ into itself and let $(m)$
be a reference measure on $X$. Consider a partition
$\{\Delta_{i}\}_{i=1}^{W}$ of the phase space $X$ and associate to
it a $W$-state Markov chain with transition probabilities
$$ p_{ij} := \frac{m(\map^{-1}\Delta_{j}\cap\Delta_{i})}{m(\Delta_{i})} .$$
It has been conjectured that a number of statistical features of
the map $\map$ can be obtained in the limit as the diameter of
the partition $\{\Delta_{i}\}_{i=1}^{W}$ vanishes of the
corresponding features of the above Markov chains (see a review
of recent results in this field and generalizations of this
procedure in \cite{Bl-mon}).

\section{Multicomponent dynamical systems and phase
transitions}\label{s:phase}

Let $\cN$ be a finite or countable collection of indices. Then by
$(X^\cN,\rho^{\cN},m^\cN)$ we denote a direct product of compact
metric spaces $(X_i,\rho_i)$ with given reference measures
$m_{i}$ on them, i.e. $X^\cN:=\otimes_{i\in\cN}X_i$,
$\rho^{\cN}(\bar x, \bar y):=\max\rho_{i}(x_{i},y_{i})$, and
$m^\cN:=\otimes_{i\in\cN}m_i$, and consider a nonsingular map
$\map^{(\cN)}:X^\cN\to X^\cN$ from this space into itself.
Consider also a collection of ``local'' maps $\map_i$ acting on
the $i$-th copy $(X_i,\rho_{i},m_i)$ of our ``local'' phase space.
We shall call the pair $(\map^{(\cN)},X^\cN)$ a {\em
multicomponent dynamical system} if its action can be decomposed
as a superposition of an ``interaction'' and the direct product
of ``local'' maps: $\cmap^{(\cN)}=\cI^{(\cN)}\circ\map^\cN$, where
$\map^\cN:=(\otimes_{i\in\cN}\map_i)$ is a direct product of maps
(from our collection of $\{\map_i\}_{i\in\cN}$). In the most
interesting spatially homogeneous case (when all spaces
$(X_{i},\rho_{i},m_{i})$ coincide) we assume also that the map
$:X^\cN\to X^\cN$ which describes the ``interaction'' between
local components (systems $(\map_i,X_i)$) of the multicomponent
dynamical system should be identical on the ``diagonal'' set
$$ X^{\cN}_{{\rm diag}} := \{\bar x \in X^{\cN}:
                             \quad x_{i}=x_{j} \; \forall i,j\}. $$
The reason of this assumption is that when all the coordinates of
a point $\bar x\in X^{\cN}$ are the same the interaction cannot
change them, i.e. $\cI^{(\cN)}\bar x=\bar x$ for any $\bar x\in
X^{\cN}_{{\rm diag}}$. For a multicomponent dynamical system with
a finite number of components this property can be considered as
a substitute for a translation invariance.

The first part of this definition, in fact, is not really
restrictive: any map $\bar\map$ from $X^\cN$ into itself can be
represented as $\bar\map\equiv
\bar\map\circ(\otimes_{i\in\cN}$Id), where Id is an identical
map. On the other hand, this representation certainly contradicts
to the second part, which assumes that the ``interaction''
($\bar\map$ in this representation) cannot change identical
elements. Note that this is one of major differences between
multicomponent models considered in statistical physics, where
the main object of interest is the evolution of systems of
particles (and therefore the ``local dynamics'' is not defined at
all, since the dynamics of individual particles without
interactions with others is trivial), and in the dynamical
systems theory.

A typical example of an admissible ``interaction'' is a space
homogeneous finite range (depending on a finite number $2K$ of
neighbors) coupling:%
\beq{adm-coupl}{
   (\cI_\ep\bar x)_i := (1-\ep)\bar x_i + \ep\sum_{j=-K}^{K}a_{j}\bar
x_{i+j} ,}%
where the parameter $\ep>0$ describes the ``interaction''
strength and $a_{i}\ge0$, $\sum_{i=-K}^{K}a_{i}=1$ are constants
defining the interaction. Observe that formula (\ref{adm-coupl})
describes a convex hall of values of coordinates of $\bar x$ in
the $K$-neighborhood of the $i$-th coordinate.

One might argue that our definition of the multicomponent system
does not cover the case when the interaction is defined in terms
of original vectors $\bar x$ rather than in terms of their images
under dynamics, i.e. $\cmap\bar x$. In particular, there is an
ubiquity of important examples \cite{discr} of multicomponent
systems obtained under space discretizations of partial
differential equations, where the system has the following form:
$$ (\bar\map_\ep \bar x)_i
 = (1-\ep)\map_i \bar x_i + \ep\sum_{j=-K}^{K}a_{j}\bar x_{i+j} ,$$
i.e. the interaction acts on $\bar x$ rather than on $\cmap\bar
x$. Formally we cannot decompose this system into the
superposition of the local dynamics described by the maps $\map_i$
and the interaction $\cI_\ep$. However we can construct an
equivalent system, which satisfies our definition, by ``doubling''
of the local systems. In the new system for each local component
$(\map_i,X_i)$ we consider its `delayed' copy acting on the phase
space $Y_i\equiv X_i$ so that the local dynamics of the i-th pair
of coordinates is defined as $\tilde\map: (\bar x_i, \bar y_i) \to
(\map_i \bar x_i, \bar x_i)$. The interaction map $\Phi(\bar x,
\bar y):=(\bar x',\bar y')$ is defined as
$$ \bar x_i' := (1-\ep)\bar x_i + \ep\sum_{j=-K}^{K}a_{j}\bar y_{i+j},\CR
   \bar y_i' := (1-\ep)\bar y_i + \ep\sum_{j=-K}^{K}a_{j}\bar y_{i+j} .$$
Thus the projection of the dynamics of the multicomponent system
$(\Phi\circ\tilde\map,~\otimes(X_i\otimes Y_i))$ to its
$x$-components coincides with the system $(\bar\map_\ep, ~\otimes
X_i)$.

An important example of formula (\ref{adm-coupl}) is the so called
diffusive coupling:
$$ (\cI_\ep\bar x)_i
 = (1-\ep)\bar x_i + \frac{\ep}3(\bar x_{i-1} + \bar x_i + \bar x_{i+1}) ,$$
i.e. the space homogeneous finite range coupling with $K=1$ and
$a_{i}\equiv1/3$. This case corresponds to the discretization of
Laplacian, indeed, we have:
$$ -\ep\bar x_i + \frac{\ep}3(\bar x_{i-1} + \bar x_i + \bar x_{i+1})
 = \frac\ep3((\bar x_{i+1} - \bar x_{i}) - (\bar x_{i} - \bar x_{i-1}))
 = \frac\ep3 (\nabla \bar x)_{i} .$$

Our definition of a natural measure is well defined for a
multicomponent dynamical system if a number of components
$|\cN|<\infty$. In this case the triple $(X^\cN, \rho^{\cN},
m^\cN)$ is again a finite dimensional metric space with a certain
reference measure and thus we can use all previous definitions.
In the case of infinite dimension ($|\cN|=\infty$) the problem is
that any probabilistic measure absolutely continuous with respect
to $m^\cN$ should coincide with it, which does not give much
freedom in the choice of initial measures. Therefore to be able
to work with infinite dimensional multicomponent dynamical
systems we need to modify the way how we choose initial measures.

Let $\cL\subseteq\cN$ be a subset of the set of indices $\cN$.
Denote by $\pi_\cL^*:\cM(X^\cN)\to\cM(X^\cL)$ -- the projection
operator in the space of probabilistic measures, defined as
$\pi_\cL^*\mu:=\int\mu~d(\otimes_{i\in(\cN\setminus\cL)}m_i)$ for
any measure $\mu\in\cM(X^\cN)$. Since there is a natural
enclosure of spaces $\cM(X^\cL)$ into the space $\cM(X^\cN)$ we
can choose a family of metrics $\dist=\dist_\cL$ acting on all
considered spaces of measures such that for any measure
$\mu\in\cM(X^\cN)$ we have
\beq{dist}{\dist(\pi_\cL^*\mu,\mu)\toas{|\cL|\to|\cN|}0.}

We shall say that a measure $\mu\in\cM(X^{\cN})$ is {\em smooth}
if its marginals $\mu_i:=\pi_{\{i\}}\mu$ are absolutely
continuous with respect to the reference measures $m_i$ for any
$i\in\cN$ and define an infinite dimensional generalization of the
natural measure $\mu_{\map}$ as a common limit of Cesaro means
$\frac1n\sum_{k=0}^{n-1}({\cmap^{(\cN)}}^*)^n\mu$ of all smooth
measures $\mu\in\cM(X^{\cN})$ having a support in a direct product
$\otimes_{i\in\cN}U_{i}$ of some open sets $U_i\subseteq X_i,
~i\in\cN$. Here $\otimes_{i\in\cN}U_{i}$ plays the role of the
basin of attraction of the measure $\mu_{\map}$.

To define the notion of phase transition we consider a family of
multicomponent systems $(\map[\gamma],X)$ depending on a certain
parameter $\gamma$. We shall say that this family has {\em the
phase transition} at the point $\gamma=\gamma_0$ if the number of
natural measures changes when the parameter $\gamma$ crosses the
value $\gamma=\gamma_0$. Observe that in the case of
multicomponent systems this may happen in two different ways.
Assume that for $\gamma<\gamma_0$ each finite dimensional
approximation $(\cmap^{(\cL)}[\gamma],X^\cL)$ have a finite
number $N[\gamma]$ of natural measures $\mu_\cL[\gamma]$ which
does not depend on $\cL$, and that any natural measure of the
complete system $\mu_\cN[\gamma]$ is a (weak) limit of measures
$\mu_\cL[\gamma]$. The main way how the phase transition may
happen is that for all $\cL$, such that $|\cL|$ is large enough,
each system $(\cmap^{(\cL)}[\gamma],X^\cL)$ goes through the
phase transition as the value of $\gamma$ crosses $\gamma_0$,
i.e. the number $\mu_\cL[\gamma]$ changes, and the same happen to
the complete system. However, in the infinite dimensional case
there is also another possibility: finite dimensional
approximations do not demonstrate any phase transition, but their
limit points either fail to correspond to the natural measures of
$\map^{(\cN)}$ at the parameter value $\gamma_0$, or a new
natural measure of the complete system appears which does not
belong to the set of limit points of natural measures of finite
dimensional approximations.

With a slight abuse of notation we denote by $\cmap^{(\cL)}$ a
$|\cL|$-dimensional approximation of our infinite-dimensional map
$\cmap^{(\cN)}$ for a given finite subset $\cL\subset\cN$ of the
set of indices. To define this approximation explicitly we need
to take into account boundary conditions, namely we have to
choose the states on the remaining infinite-dimensional part of
the phase space $X^{\cN\setminus\cL}$. Note that this can be done
in various ways. Two of them are the most common ones: fixed
boundary conditions, when the corresponding coordinates of the
vector $(\bar x)_j$ for $j\in(\cN\setminus\cL)$ are preserved at
given values, and periodic boundary conditions. The following
simple result shows that the choice of boundary conditions can
change even very rough characteristics of the dynamics.

\begin{lemma}\label{l:boundary} Let $X_i=[0,1]$ and
$\map_ix=x+\frac16x(x-1)^2$ for all $i\in\cN=\IZ^1$ and let
$\cI_\ep$ be the diffusive interaction. For a finite subset of
integers $\cL$ denote by $\cI_{\ep}^{(\cL,y)}$ the finite
dimensional approximation of the interaction with boundary
conditions fixed at the value $y$ for all coordinates not
belonging to the set $\cL$. Then for any $0<\ep<1/2$ the
multicomponent system $\cI_{\ep}^{(\cL,y)}\map^\cL,[0,1]^\cL)$
has the only one attractor (and thus the only one ergodic \SBR
measure) if $y=0$ and has $2^{|\cL|}$ attractors and ergodic \SBR
measures if the boundary condition $y=1$ and $\ep>0$ is small
enough.
\end{lemma}

\proof The map $\map_ix=x+\frac16x(x-1)^2$ has two fixed points:
the stable fixed point at $0$ and the neutrally unstable fixed
point at $1$. A straightforward calculation shows that in the
case of zero boundary conditions $y=0$ for any $\ep>0$ the fixed
point $0$ remains the only global attractor of the system, while
in the case of $y=1$ for any $\ep>0$ the fixed point at $1$
becomes stable and if $\ep<1/2$ then the fixed point at $0$
remains stable as well. \qed

The main question we shall be interested in this section is the
possibility that for each finite subset of indices $\cL$ (and a
certain choice of boundary conditions) the corresponding finite
dimensional approximation $\cmap^{(\cL)}$ has only one natural
measure, while the entire system has several of them.

In the statistical physics literature there are numerous examples
of multiparticle systems when this phenomenon takes place. One
important class of such examples is the so called voter models in
cellular automata theory. In this case each local component has
only two states and no ``local'' dynamics: the behavior of the
system is described in terms of the random ``interaction'',
namely the future state of the local coordinate is determined by
the present state of a certain number of its neighbors (including
the local component itself) with a small random error. Therefore
according to our definition this is not a multicomponent
dynamical system. The presence of phase transition for these
system has been shown by \cite{Toom} under certain assumptions on
the interaction. Later a deterministic version of this model has
been considered in \cite{MacK}. Observe that any voter model is a
probabilistic cellular automaton. Therefore by making use of the
argument from Section 3 one can immediately construct the
corresponding deterministic one-dimensional Markov map describing
this process. In fact, the construction in \cite{MacK} basically
follows this idea.

On the other hand, to the best of our knowledge, examples of
multicomponent dynamical systems with phase transitions are not
known and the only promising candidate for that is the so called
case of ``mean field'' interaction, when each local subsystem
interacts with all others (see, e.g. \cite{Kan}). In what follows
we shall give sufficient conditions under which phase transitions
cannot occur.

\begin{theorem}\label{t:stability} Assume that for any finite
subset $\cL\subset\cN$ of the set of indices there exists the
only one natural measure $\mu_\cL$ of the induced map
${\cmap^{(\cL)}}^*:=\pi_\cL {\cmap^{(\cN)}}^*\pi_\cL$, and there
exists a constant $C<\infty$ and two functions
$\phi,\psi:\IR^1\to\IR^1$ such that for any two smooth measures
$\mu,\nu\in\cM(X^\cN)$ and any two finite subsets
$\cL\subset\cL'\subset\cN$ we have
\beq{c-mult-1}{\dist({\cmap^{(\cL)}}^*\mu, {\cmap^{(\cL)}}^*\nu)
               \le C\dist(\mu,\nu) ,}
\beq{c-mult-2}{\dist({\cmap^{(\cL)}}^*\mu, \pi_\cL^*{\cmap^{(\cL')}}^*\mu)
       \le \psi(|\cL|) \toas{|\cL|\to\infty} 0 ~~\mbox{uniformly on}~ \mu ,}
\beq{c-mult-3}{\dist({\cmap^{(\cL)}}^{*^n}\mu, \mu_\cL)
       \le \phi(n) \toas{n\to\infty} 0 ~~\mbox{uniformly on}~ \mu, \cL .}
Assume also that for any measure $\mu\in\cM(X^\cN)$ we have
\beq{c-mult-4}{\dist({\cmap^{(\cL)}}^*\mu, {\cmap^{(\cN)}}^*\mu)
               \toas{|\cL|\to|\cN|}0 ,}
then the multicomponent dynamical system $\cmap^{(\cN)}$ also has
the only one natural measure.
\end{theorem}

\proof We will show that there is a weak limit of the sequence of
natural measures $\mu_*:=\mu_\cL$ as $\cL\to\cN$ and this limit
is the only one natural measure of the map $\cmap^{(\cN)}$, i.e.
$\mu_*:=\mu_{\cmap^{(\cN)}}$.

Consider a sequence of growing enclosed finite subsets
$\cL\subset\cN$. For any two finite subsets
$\cL\subset\cL'\subset\cN$ from this sequence and any positive
integer $n$ by the triangle inequality we have%
\bea{
  \dist(\mu_\cL,\pi_\cL^*\mu_{\cL'})
  \a \le \dist(\mu_\cL, {\cmap^{(\cL)}}^{*^n}\mu_{\cL'})
     + \dist({\cmap^{(\cL)}}^{*^n}\mu_{\cL'},
             \pi_\cL^*{\cmap^{(\cL')}}^{*^n}\mu_{\cL'}) \\
  \a \le \phi(n)
   + \dist({\cmap^{(\cL)}}^*\circ{\cmap^{(\cL)}}^{*^{n-1}}\mu_{\cL'},
     ~\pi_\cL^*{\cmap^{(\cL')}}^*
                \circ{\cmap^{(\cL)}}^{*^{n-1}}\mu_{\cL'}) \\
  \a \hskip 1.3cm
   + \dist(\pi_\cL^*{\cmap^{(\cL')}}^*
           \circ{\cmap^{(\cL)}}^{*^{n-1}}\mu_{\cL'},
           \pi_\cL^*{\cmap^{(\cL')}}^*
           \circ{\cmap^{(\cL')}}^{*^{n-1}}\mu_{\cL'}) \\
  \a \le \phi(n) + \psi(|\cL|)
  + C\dist({\cmap^{(\cL)}}^{*^{n-1}}\mu_{\cL'},
           \pi_\cL^*({\cmap^{(\cL')}}^{*^{n-1}}\mu_{\cL'}) \\
  \a \le \dots \le \phi(n) + \frac{C^n-1}{C-1}\psi(|\cL|) .}%
Therefore the sequence of unique natural measures $\{\mu_\cL\}$
is fundamental and thus there is a subsequence $\{\cL_i\}_i$ such
that $\mu_{\cL_i}\toas{i\to\infty}\mu_*$. By making use of the
relation~(\ref{c-mult-4}), we see that the limit measure $\mu_*$
is an invariant measure of the map $\cmap^{(\cN)}$. Thus it
remains to prove that this measure is unique and natural.

Taking a limit $\cL'=\cL_i\to\cN$ in the previous inequality we get
$$ \dist(\mu_\cL,\pi_\cL^*\mu_*)
   \le \phi(n) + \frac{C^n-1}{C-1}\psi(|\cL|) .$$
For a given smooth measure $\mu\in\cM(X^\cN)$ using the same
argument as above we get
$$ \dist({\cmap^{(\cL)}}^{*^n}\mu, \pi_\cL^*{\cmap^{(\cN)}}^{*^n}\mu)
  \le \const\phi(n) + \const^n\psi(|\cL|) ,$$
which can be made to be arbitrary small by a proper choice of
$n,|\cL|\to\infty$. Now using the assumption~(\ref{c-mult-3}) we
get that the measure $\mu_*$ is indeed the natural measure for
$\map_\cN$. \qed

\begin{corollary}\label{c:no-phase} Let a family of multicomponent
systems $\cmap[\gamma]$ satisfy the assumptions of
Theorem~\ref{t:stability} uniformly in
$\gamma\in[\gamma_1,\gamma_2]$. Then there are no phase
transitions in the interval $[\gamma_1,\gamma_2]$.
\end{corollary}

Observe that in the proof of this theorem we never used the
(spatial) decomposition of the multicomponent system into the
local components and interaction. Assume now that local
components are identical.

\begin{theorem}\label{t:stab-id} Assume that $X_i\equiv X$,
$\map_i\equiv\map$ for all $i\in\cN$, and the map $\map$ is
nonsingular w.r.t. the reference measure $m_i\equiv m$. Assume
also that the interaction $\cI$ is local, i.e. the value of $(\cI
\bar x)_i$ depends only on a finite number of `neighboring'
components $x_i$ of the vector $\bar x$. Then the statement of
Theorem~\ref{t:stability} remains valid if instead of the
assumptions (\ref{c-mult-1}, \ref{c-mult-2}) we shall use simpler
assumptions: \beq{c-mult-1'}{\dist(\cI^*\mu,\cI^*\nu) \le \const
                \dist(\mu,\nu)} \beq{c-mult-2'}{\dist({\cI^{(\cL)}}^*\mu,
                     ~\pi_\cL^*{\cI^{(\cL')}}^*\mu
                \le \psi(|\cL|)\toas{|\cL|\to\infty} 0 ,}
for any $\cL'\supset\cL$ and any $\mu,\nu\in\cM(X^\cN)$, preserve
the assumption~(\ref{c-mult-3}) and drop (\ref{c-mult-4}).
\end{theorem}

\proof Observe that ${\cmap^{(\cL)}}^*={\cI^{(\cL)}}^*{\map^\cL}^*$.
Therefore, using (\ref{c-mult-1'}), we get%
\bea{
 \dist({\cmap^{(\cL)}}^*\mu, {\cmap^{(\cL)}}^*\nu)
 \a \le \dist({\cI^{(\cL)}}^*({\map^\cL}^*\mu),
~{\cI^{(\cL)}}^*({\map^\cL}^*\nu)) \\
 \a \le \const\dist({\map^\cL}^*\mu, ~{\map^\cL}^*\nu)
 \le \const'\cdot\dist(\mu,\nu) .}%
A similar argument shows that the inequality (\ref{c-mult-2}) is
satisfied if (\ref{c-mult-2'}) holds:%
\bea{
 \dist({\cmap^{(\cL)}}^*\mu, \pi_\cL^*{\cmap^{(\cL')}}^*\mu)
 \a \le \dist({\cI^{(\cL)}}^*({\map^\cL}^*\mu),
           ~\pi_\cL^*{\cI^{(\cL')}}^*({\map^{\cL'}}^*\mu)) \\
 \a \le \psi(|\cL|) \toas{|\cL|\to\infty}0 .}%
It remains to show that the inequality~(\ref{c-mult-4}) is
satisfied automatically under our assumptions. It again follows
from the decomposition into local maps and interaction by making
use of (\ref{dist}) and (\ref{c-mult-2'}) that%
\bea{
 \dist({\cmap^{(\cL)}}^*\mu, ~{\cmap^{(\cN)}}^*\mu)
 \a = \dist({\cI^{(\cL)}}^*({\map^\cL}^*\mu),
        ~\pi_\cL^*{\cI^{(\cN)}}^*({\map^\cN}^*\mu)) \\
 \a \le \const\dist({\map^\cL}^*\mu, {\map^\cN}^*\mu)
 \toas{|\cL|\to|\cN|}0. }%
\qed

\end{document}